\documentclass[12pt]{amsart}
\usepackage[latin1]{inputenc}   
\usepackage{setspace}
\usepackage{amsthm}
\usepackage{amssymb, amsmath,psfrag}
 \usepackage[dvips]{epsfig}
\title
{An overview of harmonic analysis and the shifted wave equation on symmetric graphs}

\author{Alaa Jamal Eddine}

\address{}

\email{alaa.jamal-eddine@univ-orleans.fr}

\subjclass[2000] {Primary 44A12, 43A90;\\
\hspace{3cm}Secondary 39A12}


\newtheorem{theorem}{Theorem}[section]
\newtheorem{proposition}[theorem]{Proposition}
\newtheorem{lemma}[theorem]{Lemma}
\newtheorem{corollary}[theorem]{Corollary}

\newcommand{\bb}{\mathbb}
\newcommand{\dsps}{\displaystyle}
\newcommand\1{1\hskip-1mm\text{\rm I}}
\newcommand{\C}{\mathbb{C}}

\renewcommand\Im{\operatorname{Im}}
\newcommand{\N}{\mathbb{N}}
\renewcommand{\L}{\mathcal{L}}

\newcommand{\R}{\mathbb{R}}
\newcommand\ssb{\hskip-.25mm}
\newcommand\ssf{\hskip.25mm}
\newcommand\supp{\operatorname{supp}}
\newcommand{\Z}{\mathbb{Z}}
\makeatletter
\newenvironment{proof*}[1][\proofname]{\par
  \pushQED{\qed}%
  \normalfont \partopsep=\z@skip \topsep=\z@skip
  \trivlist
  \item[\hskip\labelsep
        \itshape
    #1\@addpunct{.}]\ignorespaces
}{%
  \popQED\endtrivlist\@endpefalse
}
\begin{document}

\begin{abstract}

Let $\Gamma$ be a  symmetric graph of type $k$ and order $r$ , where $k,r\ge 2$, and let $G$ a group acting isometrically and simlpy transitively on $\Gamma$. In this paper we give  explicite expressions of the horocyclic Abel transform and its dual as well as their inverses on $\Gamma$. We then derive the Plancherel measure for the Helgason--Fourier transform on $G$ and give a version of the Kunze-Stein phenomenon thereon. Finally, we compute the solution to the shifted wave equation on $\Gamma$, using \` Asgeirsson's mean value theorem and the inverse dual Abel transform.
\end{abstract}

\maketitle

\section{Introduction}
The Abel transform of a radial function $f$ was introuced on the hyperbolic disc $D$ by S. Helgason \cite{hel1}. It is defined as the average on the horocycles. The horocycles in the disc $D$ are the circles contained in $D$ which are tangent at the boundary $\partial D$. If $\omega\in \partial D$, the family of horocycles tangent at $\omega$ coincides with the set of level sets of the Poisson kernel $P(x,\omega)$ regarded as a function of $x$ in $D$. The Abel transform is a useful tool in harmonic analysis on the group of isometries $SL(2,\bb R)$ of the hyperbolic disc. It allows the computation of the Plancherel measure of the principal series of representations of $G$. 

The theory of spherical functions, representations and orthogonal polynomials was first generalized to the case of free groups and more generally to groups acting isometrically and simply transitively on a homogeneous tree $\bb T$ in analogy with the symmetric spaces of rank one (see \cite{ftp}, \cite{ftn}, \cite{bp}, \cite{bp1}, \cite{ca}), and then to the case of symmetric graphs (see \cite{ip}, \cite{ip1}, \cite{ks}, \cite{fp}).  More precisely in \cite{bp1}, Betori and Pagliacci defined the Abel transform on semihomogeneous trees $\bb T_{r,s}$. They gave an explicite expressions of the Abel transform and its inverse on $\bb T_{r,s}$. In the case of homogeneous trees, they proved that the Fourier transform factors into the Abel transform and the Euclidean Fourier transform as in the symmetric spaces case. Finally they gave an alternative proof of the Plancherel theorem using the Abel transform.

Our first goal is to elaborate a complete study on the Abel transform on symmetric graphs in analogy with \cite{bp1}. More precisely, let $\Gamma$ be a symmetric graph of type $k$ and order $r$, where $k,r\ge 2$. In section 3, we give  explicite expressions of the Abel transform $\mathcal{A} f$ of radial functions $f$ on $\Gamma$ and its inverse $\mathcal{A}^{-1} g$, where $g$ are even finitely supported functions on $\Z$. We then study the mapping properties of the Abel transform $\mathcal{A}$ and characterise the Abel transforms of radial finitely supported functions on $\Gamma$. We also show that the Abel transform may still be well defined on some Schwartz type function spaces.

Using the definition of the Abel transform, we derive an explicite expression of the dual Abel transform and its inverse. The dual Abel transform has been used extensively, in order to derive an explicite expression of the solution to the wave equation, (see \cite{hel} for the symmetric space case, \cite{amps} for the Damek--Ricci spaces and for the homogeneous trees cases).

A second goal of this paper is to make an overview of the harmonic analysis on the graph $\Gamma$, and give an application in this case to the results obtained later. More precisely, in section 4,  we show that the spherical Fourier transform factors as the composition of the Abel transform and the Fourier transform on $\Z$. We then give a Plancherel and inversion formulae for the Helgason Fourier transform on $\Gamma$, that is the Fourier transform for non necessarily radial functions. Recall that the Plancherel measure for the spherical Fourier transform was computed by Faraut and Picarrdello in \cite{fp} and by  Kuhn and Soardi in \cite{ks}. Our method relies on representation theoretic arguments and can be derived in a natural way from the inversion formulae for the spherical transform. We will then be able to give a version of the Kunze-Stein phenomenon for groups acting isometrically and simply transitively on $\Gamma$. 

The final goal of this paper is to solve the following shifted wave equation on $\Gamma$~:~
\begin{equation}
\begin{cases}
\beta\,\mathcal{L}^{\bb Z}_n\,u(x,n)=\,(\mathcal{L}_x^{\Gamma}-(\alpha-\beta))\,u(x,n),\\
u(x,0)=\,f(x),\;\{u(x,1)-u(x,-1)\}/2=\,g(x),
\end{cases}
\end{equation}
where $\alpha,\,\beta,\,\L^{\Z}$ and $\L^{\Gamma}$ are quantities that will be defined later. Let us first make a brief historical review on this problem. 
Helgason \cite{hel} uses \`Asgeirsson's mean value theorem (see Theorem II.5.28) to solve the wave equation on Euclidean spaces $\R^d$ (see Exercise II.F.1) and the shifted wave equation on hyperbolic spaces $\mathbb{H}^d(\R)$ (see Exercise II.F.2). This work has been extended in \cite{amps} to the case of Damek-Ricci spaces and to the case of homogeneous trees. More precisely, the method in \cite{amps} consists in using first a version of the \` Asgeirsson mean value theorem adapted to the situation and then to apply the inverse dual Abel transform. In our setting, we will first give a version of \`Asgeirsson's mean value theorem and then we use the expression of the inverse dual Abel transform obtained in section 3.

\section{Symmetric graphs}
A graph $\Gamma$ is symmetric of type $k\ge 2$ and order $r\ge 2$ if every vertex $v$ belongs exactly to $r$ polygons, with k sides each, contained in the graph, with no sides and no vertex in common except $v$, and if every nontrivial loop in $\Gamma$ runs through all the edges of at least one polygon. In other words, a symmetric graph of type $k$ and order $r$ can be thought of as a homogeneous tree of order $r$ built up with polygons with polygons with $k$ sides. Notice that, if $k=2$, $\Gamma$ is a homogeneous tree of degree $2r$. 

Different notions of distance on $\Gamma$ have been introduced in \cite{ip}. We define the distance between two vertices $v_1$ and $v_2$ as the minimal number of polygons crossed by a path connecting $v_1$ with $v_2$. Here we define the length $|v|$ of a vertex $v$, with respect to a reference vertex $o$, as the distance between $v$ and $o$. We can prove easily that the group G that acts isometrically and simply transitively on $\Gamma$ is isomorphic to the free product of $r$ copies of $\bb Z/k\bb Z$ if $k>2$, while, for $k=2$, i.e when $\Gamma$ is a homogeneous tree, $G$ is isomorphic to the free product of $t$ copies of $\bb Z$ and $s$ copies of $\bb Z/2\bb Z$, where $2t+s=r.$ 

Let us denote $\mathcal{V}$ the set of vertices of $\Gamma$. There is a natural probability  measure on $\mathcal{V}$ which  we will denote by $\mu$. Every element of $\mathcal{V}$ can be identified with an element of $G$, and every polygon of the graph corresponds, under this identification to an orbit under right translations by one of the factors $\bb Z/k\bb Z$. Under this identification, we can define the length $|g|$ of an element $g\in G$ as follows. Let $a_1,....., a_r$ denote the generators of $G$. Then, an element $g\in G$ can be written as $g=a_{i_1}^{m_1}...a_{i_n}^{m_n}$ and $|g|=n$.  Denote $S(o,n)$ the sphere of centre $o$ and radius $n$ in $G$. we have 
$$|S(o,n)|=\delta(n)=\begin{cases}
1\qquad&\text{if}\;n=0\\
r(k-1)Q^{n-1}&\text{if}\;n\ge 1
\end{cases}
$$
where $Q=(r-1)(k-1)$.

Between any two points $v_1$ and $v_2$ in $\Gamma$ such that $d(v_1,v_2)=n$, there is a unique geodesic path of the form $(P_1,...,P_n)$  where $(P_i)_{1\le i\le n}$ are polygons in $\Gamma$ and such that $v_1\in P_1\setminus P_2$ and $v_2\in P_n\setminus P_{n-1}$.

A geodesic ray $\mathcal{P}$ in $\Gamma$ is a one sided sequence $\{P_n,\,n\in\bb N\}$ of polygons. We say that $v\in\mathcal{P}$ if $v\in P_n$ for some $n\in \bb N$. We denote by $\Omega$ the set of all geodesic rays starting with $o$. $\Omega$ is called the boundary of $\Gamma$. $G$ acts by left translation on words in $\Omega$. For $x\in G$, let $E(x)$ denotes the subset of $\Omega$ of words that begin with the reduced word $x$. Then, $\{E(x),\,x\in G\}$ is a base of the topology of $\Omega$ making $\Omega$ a compact topological space. Let $\nu$ denotes a probability measure on $\Omega$ defined by
$$\nu(E(x))=\frac1{\delta(n)}\;\text{if}\;|x|=n.$$
We then denote $(\Omega,\nu)$ the Poisson boundary of $G$ with respect to $\mu_1$. The group $G$ acts on measures on $\Omega$, particularily on $\nu$ by
$$\nu_x(A)=\nu(x^{- 1}A)$$
for all $x\in G$, and for all Borelian set $A\in \Omega$.  The measure $\nu_x$ is absolutely continuous with respect to $\nu$, and we have, thus for $\omega\in \Omega$, such that $\omega=a_{i_1}^{n_1}a_{i_2}^{n_2}......$ we denote $\omega_m=a_{i_1}^{n_1}...a_{i_m}^{n_m}$. The Radon Nikodym derivative
\begin{equation}\label{poisson}
\nu(x^{-1}E(\omega_m))/\nu(E(\omega_m))=Q^{\zeta(x,\omega)}
\end{equation}
where $\zeta$ denotes the Busemann function on $\Gamma$
\begin{equation*}
\zeta(x,\omega)= d(o,\omega_m)-d(x,\omega_m).
\end{equation*}
We remark that $\zeta(x,\omega)$ doesn't depend on $\omega$ if $m>|x|$. The quantity obtained in $(\ref{poisson})$ denotes the Poisson kernel $P(x,\omega)$, it verifies the following cocycle identities~:
\begin{equation}\label{cocycle}
P(o,\omega)=1\qquad\text{and}\;P(xy,\omega)=P(y,x^{-1}\omega)\,P(x,\omega).
\end{equation}
\indent
We shall need some function spaces and related notations. Given a discrete space $X$, we denote by $\mathcal{D}(X)$ the space of all finitely supported functions on $X$.  A function $f$ on $G$ is radial if it is constant on $S(o,n)$ for all $n\in\bb N$. If $E(G)$ is a space of functions on $G$, then $E(G)^{\sharp}$ will denote the subset of $E(G)$ of radial functions. If $E(\bb Z)$ is a space of functions on $\bb Z$, then $E(\bb Z)_{\text{even}}$ will denote the subspace of $E(\bb Z)$ of even functions therein.

We define the convolution product on $G$ by 
$$
f*g(x)=\sum_{y\in G}f(y) g(y^{-1}x).
$$
If $g$ is radial, we have
$$f*g(x)=\sum_{n\in\bb N}g(n)\sum_{d(x,y)=n}f(y)
$$
We use the variable constants convention, and denote by $C$ a constant who will depend only on unvariable data.
\vspace{-0.3cm}
\section{The Abel transform}
\begin{figure}[!h]
\psfrag{0}[c]{$0$}
\psfrag{1}[c]{$1$}
\psfrag{2}[c]{$2$}
\psfrag{-1}[c]{$-1$}
\psfrag{h}[c]{$h$}
\psfrag{omega}[c]{$\omega$}
\psfrag{origin}[c]{$0$}
\includegraphics[width=13cm,height=8.0cm]{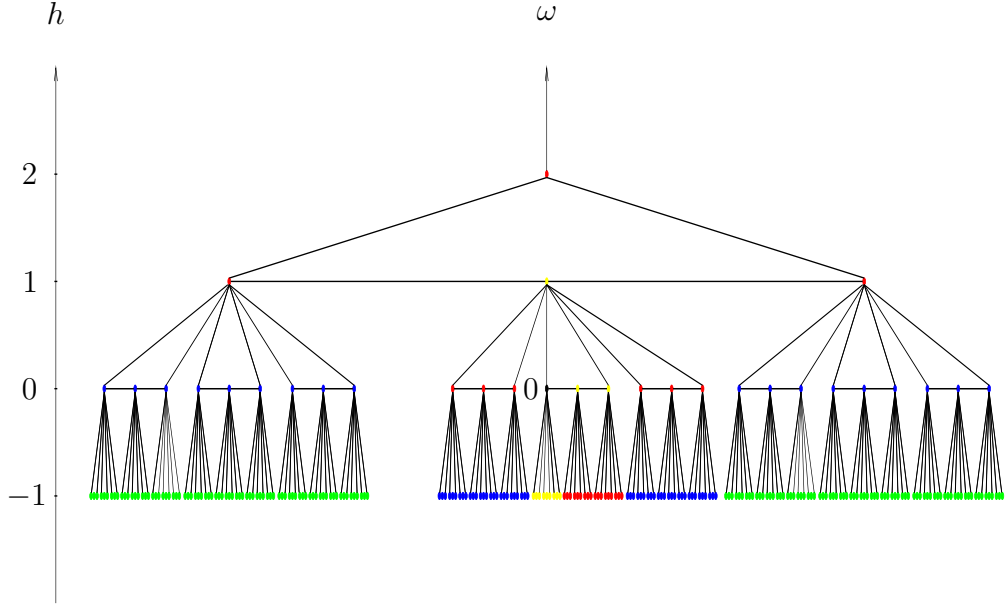}
\vspace{3mm}
\caption{Spheres in the upper-half space of  \ssf
$(\mathbb{Z}/4\ssf\mathbb{Z})^{\star\ssf4}$}\label{spheres}
\end{figure}

In this section, we give an explicite expression of the Abel transform, which is the horocyclic Radon transform on radial functions. We first begin  by describing the horocycles on $\Gamma$. Recall that $\mathcal{V}$ denotes the set of vertices of $\Gamma$ and $\Omega$ its boundary. For $\omega\in \Omega$ and $x\in \mathcal{V}$, there exists a unique geodesic ray issued from $x$ and joining $\omega$ that we will denote by $[x,\omega]$. Let $x,y\in \Gamma$ and $\omega\in \Omega$. We denote by $z$ the confluence point of the geodesic rays $[x,\omega],[y,\omega]$, that is the last point on $[x,\omega]$ laying on the geodesic ray $[y,\omega]$. We then define $\zeta_{\omega}(x,y)=d(x,z)-d(y,z)$. Thus, for $\omega$ in $\Omega$, the relation $\zeta_{\omega}(x,y)=0$ defines an equivalence relation on $\mathcal{V}$,  and the equivalence classes denote the horocycles of $\Gamma$. 
The choice of $o$ as an origin will permit us to enumerate the horocycles and to prove that the set of horocycles $H_h(\omega)$ is one-to-one with $\Omega\times\bb Z$. More precisely, for $\omega\in\Omega$ and $h\in\Z$, we have 
$$H_h(\omega)=\{ x\in\mathcal{V},\,\zeta_{\omega}(o,x)=h\}.$$
The function $\zeta_{\omega}(o,x)$ is nothing else than the Busemann function $\zeta(x,\omega)$ defined above. Thus, horocycles are the level sets of the Poisson kernel. We may prove easily that $\mathcal{V}$ decomposes disjointly as~:
$$\mathcal{V}=\dsps\mathop{\bigcup}_{h\in\bb Z}H_h(\omega),$$
in particular, $o\in\,H_0(\omega)$.
The horocyclic Radon transform is defined on $\mathcal{D}(\mathcal V)$ by
$$
\mathcal{R}f(\omega,h)=\dsps\sum_{x\in H_h(\omega)}f(x)
$$
and the Abel transform on $\mathcal{D}(\mathcal{V})^{\sharp}$ by
$$
\mathcal{A}f(\omega,h)=Q^{\frac{h}2}\,\mathcal{R}f(\omega,h).
$$
The following proposition gives an explicit formula for the Abel transform of a radial function $f$. Such a function may be identified with a function on $\bb N$, and we denote by $f(n)$ the common value $f(x)$ when $|x|=n$. 
\begin{proposition}\label{abelprop}
Denote $\sigma=k-2$.
If $f\in\mathcal{D}(\mathcal V)^{\sharp}$, then
$$\mathcal{A}f(\omega,h)= Q^{\frac{|h|}{2}}\,f(|h|)\,+\,\sigma\dsps\sum_{j\geq 1}Q^{\frac{|h|}{2}+j-1} f(|h|+2j-1)\,+\,\frac{r-2}{r-1}\,\dsps\sum_{j\geq 1}Q^{\frac{|h|}{2}+j}\,f(|h|+2j)$$
$\forall\,(\omega,h)\in\Omega\times\bb Z$. 
\end{proposition}
\noindent
Consequently, $\mathcal{A}f$ is constant in the first variable, and is even in the second.  The following lemma is crucial in the proof of this proposition.
\begin{lemma}
Given $\omega\in\Omega$. For all $n\in\bb N$ and $h\in\bb Z$, let $b(n,h)=\operatorname{Card}\left\{H_h(\omega)\bigcap S(o,n)\right\}.$  Then, 
$$b(n,h)=
\begin{cases}
0 \qquad&\textrm{if}\;n< |h|\\
Q^{-h_{-}}&\textrm{if}\;n=|h|
\\\sigma\,Q^{-h_{-}+j-1}&\textrm{if}\;n=|h|\!+\!2j\!-\!1\;\textrm{where}\;j\ge 1
\\(r-2)(k-1)\,Q^{-h_{-}+j-1}&\textrm{if}\;n=|h|\!+\!2j\;\textrm{where}\;j\ge 1,
\end{cases}
$$
where $h_{-}=min(0,h)$.
\end{lemma}
\begin{proof*}[Proof of lemma]
The proof is similar to that on semihomogeneous trees (see \cite{bp1}).  Given $\omega\in\Omega,$ and let $\omega=a^{n_1}_{i_1}a^{n_2}_{i_2}.....$ be the reduced word representation of $\omega$. We denote by $\omega_0=a^{n_1}_{i_1},$ $\omega_1=\omega_0\, a_{i_2}^{n_2},...$ the nodes of $[o,\omega]$, that is the points on $[o,\omega]$ where we change of horocycles. If $h\ge 0$, $\omega_h$ is the only element of $H_h(\omega)\bigcap S(o,h)$, thus $b(h,h)=1.$ Now, if $h<0$, a point $x$ belongs to $H_h(\omega)\bigcap S(o,|h|)$ if and only if the origin $o$ is the confluence point of the geodesics $[x,\omega]$ and $[o,\omega]$, or
$$
\operatorname{Card}\{x\in S(o,-h)|\,o\notin [x,\omega]\cap [o,\omega]\}=(k-1)^{-h}\,(r-1)^{-h-1},
$$
then $b(-h,h)=\operatorname{Card}\{S(o,-h)\}-(k-1)^{-h}\,(r-1)^{-h-1}=Q^{-h}.$ Clearly $b(n,h)=0$ when $|h|>n$. In the other cases, we proceed by reccurence on $n$, using the following equality
\begin{equation}\label{bnh}
b(n+2,h)=Q b(n,h)
\end{equation}
in both cases. We consider first the case $h=n-1$. In this case $h$ and $n$ are of different parity. The intersection points of $S(o,n)$ and $H_h(\omega)$ are situated on the polygone that has $\omega_h$ and $\omega_{h+1}$ as nodes, thus $\omega_h\in H_h(\omega)\setminus S(o,n)$, and $\omega_{h+1}\in H_{h+1}(\omega)$. Clearly the other points of this polygon belongs to $S(o,n)\bigcap H_h(\omega)$, so  $b(n,n-1)=\sigma,$ thus we have proven the formula for $h=n-1$. If $h=n-2$, $h$ and $n$ have the same parity. In this case, we consider the points $v\in\mathcal{V}$ such that $|v|=n=|\omega_h+2|$, and $v\in H_h(\omega)$. These points are vertices of polygons issued from $\omega_{h+1}$. Or, all vertice of polygons situated between $\omega_{h+1}$ and $\omega_{h+2}$ belong to $H_{h+1}(\omega)$. It's clear that the $(r-2)$ other polygons have vertices that belong to $S(o,n)\bigcap H_h(\omega)$ for $h=n-2$ (see fig $\ref{spheres}$). Then, if $h=n-2$, $b(n,h)=(r-2)(k-1).$ We conclude using $(\ref{bnh})$.
\end{proof*}

 By definition, the Abel transform of a radial function $f$ may be written as
$$
\mathcal{A}f(\omega,h)=Q^{\frac{h}2}\sum_{n\in\bb N}b(n,h)\,f(n),
$$
thus we prove the proposition $\ref{abelprop}$ by substituing the value of $b(n,h)$ in this formula.

\begin{proposition}
For all $k\le r$, the Abel transform $\mathcal{A}:\mathcal{D}(G)^{\sharp}\to \mathcal{D}(\bb Z)_{even}$ is an isomorphism, and its inverse is given by
\begin{align}\label{inv}
\mathcal{A}^{-1}g(n)=\frac{1}{k}\,Q^{-\frac{n-1}{2}}&\Big\{\sum_{m=1}^{\infty}Q^{-\frac{m}{2}}\,[g(n+m-1)-g(n+m+1)]\\&+\sum^{\infty}_{m=1}(-1)^{m-1}(k-1)^m Q^{-\frac{m}{2}}[g(n+m-1)-g(n+m+1)]\nonumber\Big\}.
\end{align}
which is equivalent to the following expression
\begin{align}\label{inv1}
\mathcal{A}^{-1}g(n)=\,Q^{-\frac{n}{2}}&\Big\{g(n)-(k-2)\,Q^{-\frac12}\,g(n+1)\\&-\frac{Q-1}{k}\,\sum_{m=1}^{\infty}Q^{-\frac{m}{2}}\,g(n+m)\nonumber\\&-\frac{r-k}{r}\sum_{m=1}^{\infty}(-1)^m\,(k-1)^m\,Q^{-\frac{m}{2}}\,g(n+m)\nonumber\Big\}.
\end{align}
Note that these sums are finite as long as we consider functions $g$ of finite supports.
\end{proposition}
\begin{proof*}
By Proposition $\ref{abelprop}$, we have
\begin{align*}
&Q^{-\frac{n}{2}}\Big\{\mathcal{A}f(n)-\,\mathcal{A}f(n\!+\!2)\,\!+\!\,Q^{-\frac{1}{2}}[\mathcal{A}f(n\!+\!1)\,-\,\mathcal{A}f(n\!+\!3)]\Big\}=\\&f(n)-\,f(n\!+\!2)\!+\!\,(k-1)\,[\,f(n\!+\!1)-f(n\!+\!3)\,],
\end{align*}
or $f$ is of finite support, so 
\begin{align*}
f(n)\!+\!\,(k-1)\,f(n\!+\!1)=&\,\dsps\sum_{j\geq 0} \left\{f(n\!+\!2j)-\,f(n\!+\!2j\!+\!2)\right\}\!+\!\\&(k-1)\,\dsps\sum_{j\ge 1}\left\{\,f(n\!+\!2j-1)-\,f(n\!+\!2j\!+\!1)\right\}
\end{align*}
thus
 \begin{align}\label{fn}
f(n)\!+\!\,(k-1)\,f(n\!+\!1)=&\,\dsps\sum_{j\geq 0}Q^{-\frac{n}{2}-j}\Big\{\mathcal{A}f(n\!+\!2j)-\,\mathcal{A}f(n\!+\!2j\!+\!2)\Big\}\\&\!+\! \,\dsps\sum_{j\ge 1} Q^{-\frac{n-1}{2}-j}\Big\{\mathcal{A}f(n\!+\!2j\!-\!1)\,-\,\mathcal{A}f(n\!+\!2j\!+\!1)]\Big\}\nonumber.
\end{align}
Let $F(n)=(k-1)^n f(n)$, so when multiplying both members of $(\ref{fn})$ by $(k-1)^n$, we obtain the following reccurence formula 
\begin{equation}\label{F}
F(n)\!+\!\,F(n\!+\!1)=\,(k-1)^n G(n),
\end{equation}
where
\begin{align}\label{G(n)}
G(n)&=\sum_{j\ge 0}Q^{-\frac{n}{2}-j}\{g(n+2j)-g(n+2j+2)\}\\&+\sum_{j\ge 0}Q^{-\frac{n+1}{2}-j}\{g(n+2j+1)-g(n+2j+3)\}.\nonumber
\end{align}
We rewrite $(\ref{F})$ as
$$
F(n)=(k-1)^n G(n)-F(n+1).
$$
Then
\begin{equation}\label{hn}
F(n)=\,\dsps\sum_{\ell\ge 0}(-1)^{\ell} (k-1)^{n+\ell}\,G(n\!+\!\ell).
\end{equation}
Thus
\begin{equation}\label{fnn}
f(n)=\sum_{\ell\ge 0}(-1)^{\ell} (k-1)^{\ell}\,G(n\!+\!\ell).
\end{equation}
Note that this sum is finite for functions $g=\mathcal{A}f$, consequently $G$ is of finite support. Substituing $(\ref{G(n)})$ in $(\ref{fnn})$, then leting $m=\ell+2j+1$, respectively $m=\ell+2j$, and 
\vspace{-0.3cm}
distinguishing the case $m$ even and $m$ odd, we have
\begin{align}\label{123}
\mathcal{A}^{-1}g(n)&=\,Q^{-\frac{n}{2}}\,\{g(n)-g(n+2)\}\\&-\sum_{m\,\textup{even}\,\ge 2}\frac{(k-1)^m -1}{k}\,Q^{-\frac{n+m-1}{2}}\,\{g(n+m-1)-g(n+m+1)\}\nonumber\\&+\sum_{m\,\textup{odd}\,\ge 3}\frac{(k-1)^m +1}{k}\,Q^{-\frac{n+m-1}{2}}\{g(n+m-1)-g(n+m+1)\}\nonumber\\&=\frac{1}{k}\,\Big\{\sum_{m\,\textup{odd}\,> 0}[(k-1)^m+1]\,Q^{-\frac{n+m-1}{2}}\,\{g(n+m-1)-g(n+m+1)\}\nonumber\\&-\sum_{m\,\textup{even}\,> 0}[(k-1)^m-1]\,Q^{-\frac{n+m-1}{2}}\,\{g(n+m-1)-g(n+m+1)\}\Big\}.\nonumber
\end{align}
Rearranging the terms, we obtain $(\ref{inv}),$ then $(\ref{inv1})$.
\end{proof*}
\begin{corollary}
Let $f$ be a radial function on $\Gamma$. We have the following~:~
$$\supp f\subset B^{\prime}(0,n)\;\textrm{if and only if}\;\operatorname{supp}\,\mathcal{A}f\subset\,[-n,+n].
$$
\end{corollary}
Our next goal is to extend the definition of the Abel transform to Schwartz spaces. For $0<p<\infty$, we denote by $\mathcal{S}_p(G)$ the set of functions $f:G\to\bb C$ such that, for all $m\in\N$,
\begin{equation}\label{fpm}
\left\|f\right\|_{(p,m)}=\sup_{{x\in G}}\,(1\!+\!|x|)^m Q^{|x|/p}|f(x)|<\infty.
\end{equation}
It's a Frechet space for the increasing norms family $(\ref{fpm})$. We will also consider the space $\mathcal{S}(\bb Z)_{\text{even}}$ of even functions $g:\bb Z\to\bb C$ such that, for all $m\in\N$,
$$\left\|g\right\|_{(m)}=\sup_{{n\in\bb Z}}\,(1\!+\!|n|)^m|g(n)|<\infty.
$$
\begin{proposition}
\begin{enumerate}
\item
For all $1\le p\le 2,$ the Abel transform extends to a continuous homomorphism  from $\mathcal{S}_p(G)^{\sharp}$ to $Q^{-(\frac1{p}-\frac12)|.|}\,\mathcal{S}(\bb Z)_{\text{even}}$
\item
If $k\le r$, then $\mathcal{A}^{-1}$ is a continuous homomorphism from $Q^{-(\frac1{p}-\frac12)|.|}\mathcal{S}(\bb Z)_{\text{even}}$ to $\mathcal{S}_p(G)^{\sharp}$, for all $1\le p\le 2$
\item
If $k>r$, then $\mathcal{A}^{-1}$ is a continuous homomorphism from $Q^{-\frac12|.|}\mathcal{S}(\bb Z)_{\text{even}}$ to $\mathcal{S}_1(G)^{\sharp}.$
\end{enumerate}
\end{proposition}
\begin{proof*}
Let $1\le p\le 2$. We will show that, for all $m\in \bb N$ there exist $C>0$ such that, for $f\in \mathcal{S}_p(G)^{\sharp},$
\begin{equation}\label{aa}
\left\|Q^{(\frac{1}{p}-\frac12)|.|}\mathcal{A}f\right\|_{(m)}\le C\left\|f\right\|_{(p,m+2)}.
\end{equation}
Or by hypothesis we have
$$|f(n)|\le C^{\prime}(1\!+\!n)^{-m-2}Q^{-n/p}\left\|f\right\|_{(p,m+2)}$$
for all $f\in S_p(G)^{\sharp}$ and $n\in\bb N$. Then
\begin{align*}
(1\!+\!h)^m Q^{(\frac1{p}-\frac12)h}|\mathcal{A}f(h)|&\le\,(1\,\!+\!\,h)^m\sum_{j\in\bb N} Q^{\frac{h}{p} +\frac{j}{2}}|f(h\,\!+\!\,j)|\\& \le C\sum_{j\in\bb N}(1\,\!+\!\,j)^{-2}\left\|f\right\|_{(p,m+2)}
\end{align*}
for all $h\in \bb N$, and this ends the proof of the first part.
For the second part, we will show that, for all $m\in \bb N$, there exist $C>0$ such that, for all $g\in \,Q^{-(\frac{1}{p}-\frac12)\,|.|}\mathcal{S}(\bb Z)_{\text{even}}$,
$$\left\|\mathcal{A}^{-1}g\right\|_{(p,m)}\le C\left\|g\right\|_{(m\!+\!2)}.
$$
By hypothesis, for $g\in Q^{-(\frac1{p}-\frac12)|.|}\, \mathcal{S}(\bb Z)_{\text{even}}$ and $m\in \bb N$, we have 
$$
|g(n)|\le\,(1\!+\!n)^{-m-2}\,Q^{-(\frac{1}{p}-\frac12)\,n}\,\left\|Q^{-(\frac{1}{p}-\frac12)\,|.|}g\right\|_{(m+2)}
$$
for all $n\in\N$.
Or, for $n\in\N$, we have $|\mathcal{A}^{-1}g(n)|\le\,\dsps\sum_{i\ge 0}(k-1)^i\,|G(n\!+\!i)|,$ where 
$(k-1)^i\le Q^{\frac{i}{2}}\le Q^{\frac{i}{p}}$, and
$$
|G(n\!+\!i)|\le C\dsps\sum_{j\ge 0} Q^{-\frac{n+i+j}{2}}\,|g(n\!+\!i\!+\!j)|,$$
so we deduce the following estimations~:~
\begin{align*}
(1\!+\!n\!+\!i)^{m+2}\,Q^{\frac{n+i}{p}}\,|G(n\!+\!i)|&\le C\dsps\sum_{j\ge 0} Q^{-\frac{j}{p}}\,(1\!+\!n\!+\!i\!+\!j)^{m+2}\,Q^{(\frac{1}{p}-\frac12)(n+i+j)}\,|g(n\!+\!i\!+\!j)|\\&\le C\left\|Q^{(\frac{1}{p}-\frac12)|.|}g\right\|_{(m\!+\!2)}\;\;\forall n,i,\,\in\bb N.
\end{align*}
Thus, for $n\in\N$, we have
\begin{align*}
(1\!+\!n)^m\,Q^{\frac{n}{p}}\,|\mathcal{A}^{-1}g(n)|&\le \dsps\sum_{i\ge 0} (1\!+\!i)^{-2}(1\!+\!n\!+\!i)^{m+2}\,Q^{\frac{n+i}{p}}\,|G(n\!+\!i)|\\&\le C\left\|Q^{(\frac{1}{p}-\frac12)\,|.|}\,g\right\|_{(m+2)}.
\end{align*}
We thus proved the second part. The last part can be proven in the same way.
\end{proof*}
\begin{corollary}
\begin{enumerate}
\item
If $k\le r$ , then , for all $p\in\,[1,2],$ the Abel transform is a toplogic isomorphism between $\mathcal{S}_p(G)^{\sharp}$ and $Q^{-(\frac1{p}-\frac12)|.|}\mathcal{S}(\bb Z)_{\text{even}}.$
\item
If $k>r$, then the Abel transform is a topologic isomorphism between $\mathcal{S}_1(G)^{\sharp}$ and $Q^{-\frac12|.|}\mathcal{S}(\bb Z)_{\text{even}}$.
\end{enumerate}
\end{corollary}
Next, we will derive an explicite expression of the dual Abel transform and its inverse.
The dual Abel transform $\mathcal{A}^*g$ of an even function $g:\bb Z\to\bb C$ is defined as follows~:
\begin{equation}\label{abeldef}
\sum_{n\in \bb N}\mathcal{A}^*g(n)\,f(n)\, \delta(n)=\,\sum_{h\in\bb Z} g(|h|)\,\mathcal{A}f(|h|)
\end{equation}
for all $f\in\mathcal{D}(\mathcal{V})^{\sharp}.$ Recall that $\delta(n)$ denotes the cardinal of the sphere of radius $n$, and $\sigma=k-2$.
\begin{theorem}\label{abeldu}
Let $g~:~\Z\to\C$ be an even function. Then,
 $$\mathcal{A}^*g(0)=\,g(0)$$
and, for all $n> 0$, 
\begin{equation}
\mathcal{A}^*g(n)=~\,2\,\frac{r-1}{r}\,Q^{-\frac{n}{2}}g(n)\,+\,\sigma\, \frac{r-1}{r}\,Q^{-\frac{n-1}{2}}\hspace{-0.7cm}\sum_{
-n< j< n \atop{\textrm{$j$ and $n$ of different parity}}}\hspace{-1.0cm}g(\pm j)\,+\,\frac{r-2}{r}\,Q^{-\frac{n}{2}}\hspace{-0.4cm}\sum_{-n< j< n \atop{\textrm{$j$ and $n$ of the same parity}}} \hspace{-1.0cm}g(\pm j)
\end{equation}
\end{theorem}
\begin{proof*}
\begin{equation}\label{gl}
\sum_{h\in\bb Z}\mathcal{A}f(h)\,g(h)=\,\mathcal{A}f(0)g(0)\,\!+\!\,2\,\dsps\sum_{h\ge 1}\mathcal{A}f(h)\,g(h)=\dsps\sum_{n\geq 0}Q^{\frac{n}{2}}f(n) G(n).
\end{equation}
where\, $G(0)=g(0)$\, and
$$G(n)= 2\,g(n)\,\!+\!\, \frac{r-2}{r-1}\hspace{-0.3cm}\dsps\sum_{-n< j< n\atop{j\,\text{et}\,n\,\text{of the same parity}}}\hspace{-0.5cm}g(\pm j) \,\!+\!\,\sigma\, Q^{-\frac12}\hspace{-0.4cm}\dsps\sum_{-n< j< n \atop{j\,\text{et\,$n$\,of different parity}}}\hspace{-0.7cm}g(\pm j),$$
we then conclude easily.
\end{proof*}
We will next establish the expression of the inverse dual Abel transform.
\begin{theorem}
The dual Abel transform is an isomorphism between the space of even functions $g:\bb Z\to\bb C$  and the space of radial functions $f:\mathcal{V}\to\bb C$. Its inverse is given by 
$$g(0)=\,f(0),\;\;\;g(1)=\,-\frac{\sigma}{2}\,Q^{-\frac12}\,f(0)\,\!+\!\,\frac{r(k-1)}{2}\,Q^{-\frac12}\,f(1)$$
and
\begin{align*}
g(n)=&-\frac{1}{2k}\left\{q-1\!+\!(r-k)(1-k)^n\right\}\,Q^{-\frac{n}{2}}\,f(0)\\&-\frac{r(k-1)}{2k}\dsps\sum_{0< j< n-1}\left\{q-1\!+\!(r-k)(1-k)^{n-j}\right\}Q^{j-\frac{n}2-1}\,f(j)\\&-\frac12 \,r\,(k-1)\,\sigma \,Q^{\frac{n}{2}-2}\,f(n-1)+\,\frac12\, r\,(k-1)\,Q^{\frac{n}{2}-1}\,f(n)
\end{align*}
\noindent
for all $n\ge 2$, with the usual convention that sum on empty sets is equal to zero.
\end{theorem}
\begin{proof*}
Given $g\in \mathcal{D}(\bb Z)_{\text{even}},$ let 
\begin{align*}
&G(n)=Q^{\frac{n}{2}}\,g(n)\;\;\text{and}\\&F(n)=\frac12\,r\,(k-1)\,Q^{\frac{n}{2}}\left\{Q^{\frac{n+2}{2}}\,\mathcal{A}^*f(n\!+\!2)-\,Q^{\frac{n}{2}}\,\mathcal{A}^*g(n)\right\}.
\end{align*}
From theorem $\ref{abeldu}$ we have,
\begin{equation}\label{identitG}
G(n\!+\!2)\,\!+\!\,\sigma\,G(n\!+\!1)-\,(k-1)\,G(n)=\,F(n)
\end{equation}
for all $n\in\bb N$. The equation $(\ref{identitG})$ is an inhomogeneous linear reccurence relation of second order. Its homogeneous solution is given by
$$G_{\text{hom}}(n)=\,c_1 \!+\!c_2 (1-k)^n.$$
We then use the constant variation method  to resolve the inhomogeneous equation~:
\begin{equation}\label{varcons}
G(n)=\,c_1(n)\!+\!\,c_2(n)\,(1-k)^n.
\end{equation}
with the additionnal condition~:
\begin{equation}\label{condaddi}
c_1(n\!+\!1)-\,c_1(n)=\,\left\{c_2(n\!+\!1)-\,c_2(n)\right\}\,(1-k)^{n+1}=\,0.
\end{equation}
Thus
\begin{equation}\label{Gn1}
G(n\!+\!1)=\,c_1(n)\!+\!\,c_2(n)\,(1-k)^{n+1},
\end{equation}
hence
\begin{equation}\label{Gn2}
G(n\!+\!2)=\,c_1(n\!+\!1)\!+\!\,c_2(n\!+\!1)\,(1-k)^{n+2}.
\end{equation}
Substituing $(\ref{varcons})$, $(\ref{Gn1})$, $(\ref{Gn2})$ in $(\ref{identitG})$, we obtain
\begin{equation}\label{cF}
c_1(n\!+\!1)-\,c_1(n)\!+\!\,\left\{c_2(n\!+\!1)-\,c_2(n)\right\}\,(1-k)^{n+2}=\,F(n).
\end{equation}
The system of equations formed by $(\ref{condaddi})$ and $(\ref{cF})$ have the following solution
$$
\begin{cases}
c_1(n\!+\!1)-\,c_1(n)=\,k^{-1}\,F(n)\\
c_2(n\!+\!1)-\,c_2(n)=\,-k^{-1}\,(1-k)^{-n-1}\,F(n)
\end{cases}
$$
hence
$$
\begin{cases}
c_1(n)=\,c_1\!+\!\,k^{-1}\dsps\sum_{0\le j< n} F(j)\\
c_2(n)=\,c_2-\,k^{-1}\dsps\sum_{0\le j< n}(1-k)^{-j-1}\,F(j)
\end{cases}
$$
and
\begin{equation}\label{G}
G(n)=\,c_1\!+\!\,c_2(1-k)^n\!+\!\,\dsps\sum_{0\le j\le n-2}\frac{1-\,(1-k)^{n-j-1}}{k}\,F(j).
\end{equation}
We compute the constants
$$
\begin{cases}
c_1=\,\frac12\,\mathcal{A}^*g(0)+\,\frac{r(k-1)}{2k}\,\mathcal{A}^*g(1),\\
c_2=\,\frac12\,\mathcal{A}^*g(0)-\,\frac{r(k-1)}{2k}\,\mathcal{A}^*g(1),
\end{cases}
$$
using the initial values
\begin{equation*}
\begin{cases}
G(0)=\,g(0)=\,\mathcal{A}^*g(0),\\
G(1)=\,Q^{\frac12}\,g(1)=\,\frac{r(k-1)}{2}\,\mathcal{A}^*g(1)-\,\frac{\sigma}{2}\,\mathcal{A}^* g(0).
\end{cases}
\end{equation*}
We conclude by substituting in  $(\ref{G})$ the expressions of $c_1,\,c_2$ and $F$.
\end{proof*}
We will now describe some consequences of the results obtained above. Let us recall first the definition of the spherical Fourier transform on $\Gamma$. There is a natural Laplace operator $\mathcal{L}$ on $\Gamma$, defined by 
$$\mathcal{L}f(x)=f(x)-\frac1{r(k-1)}\sum_{y:d(x,y)=1}f(y),$$
and the spherical functions are the radial eigenfunctions $\varphi$ of $\mathcal{L},$ normalized with the condition $\varphi(o)=1$. Each of these functions may be represented (see e.g., \cite{ip1}) in terms of the Poisson kernel by the formula
\begin{equation}\label{phi}
\varphi_{\lambda}(x)=\int_{\Omega}P(x,\omega)^{\frac12+i\lambda}d\omega,
\end{equation}
so that $\varphi_{\lambda}$ is an eigenfunction of the Laplacian $\L$ with eigenvalue 
$$
\gamma(\lambda)= \frac{Q^{\frac12+i\lambda}+Q^{\frac12-i\lambda}+\sigma}{r(k-1)}.
$$
The spherical Fourier transform of a function $f\in\mathcal{D}(\Gamma)^{\sharp}$  is defined by
$$
\mathcal{H}f(\lambda)=\sum_{x\in\Gamma}f(x)\varphi_{\lambda}(x)\qquad\forall\,\lambda\in\bb C.
$$
This definition of the spherical Fourier transform may be extended to some classes of functions of infinite support. Indeed, using the definition $(\ref{phi})$ of the spherical functions, we have $|\varphi_{\lambda}(x)|\le 1$ for all $x\in G$ and $\lambda\in\C$ such that $|\Im \lambda|\le \frac12.$ Hence, $\mathcal{H}f(\lambda)$ may be defined for function radial functions $f\in L^1(G)$ and for $\lambda\in\C$ such that $|\Im \lambda|\le \frac12.$

Let $\tau=\frac{2\pi}{\ln Q}.$ $\mathcal{H}f$ is even and  $\tau$--periodic. More generally, we can define the Helgason-Fourier transform of  a function $f\in\,\mathcal{D}(G)$ (not necessarily radial) by
\begin{equation}\label{hel}
\hat{f}(\lambda,\omega)=\sum_{x\in G}f(x)\,P(x,\omega)^{\frac12+i\lambda},\quad\forall\,(\omega,\lambda)\in\Omega\times\bb C.
\end{equation}
Using the expression of $P(x,\omega)$, we obtain
\begin{equation}\label{relationAF}
\hat{f}(\lambda,\omega)=\sum_{h\in\bb Z}Q^{(\frac12 +i\lambda)h}\sum_{x\in H_h(\omega)}f(x)=\,\mathcal{F}_h\left\{Q^{\frac{h}{2}}\mathcal{R}f(\omega,h)\right\}(\lambda),
\end{equation}
where $\mathcal{F}$ is the Fourier transform on $\bb Z$, given by ~:
$$\mathcal{F}g(\lambda)=\sum_{n\in\bb Z} Q^{in\lambda}g(n),$$
and its inverse is given by 
$$
g(n)=\,\frac{1}{\tau}\int_{-\frac{\tau}{2}}^{\frac{\tau}{2}}\mathcal{F}g(\lambda)\,Q^{-in\lambda}\,d\lambda.
$$
Recall that for radial functions $f$, $\mathcal{A}f$ doesn't depend on $\omega$, whence, using $(\ref{relationAF})$, we observe that, for radial functions $f$,  $\hat{f}$ doesn't depend on $\omega$.
In this case, using the definition $(\ref{phi})$ of spherical functions, we have
$$\hat{f}(\lambda,\omega)=\int_{\Omega}\hat{f}(\lambda,\omega)\,d\nu(\omega)\,=\,\sum_{x\in G}f(x)\int_{\Omega}P^{\frac12 + i\lambda}(x,\omega)\,d\nu(\omega)=\,\mathcal{H}f(\lambda).$$
This shows that, the spherical Fourier transform and the Helgason-Fourier transform coincide on radial functions, and we have
\begin{equation}\label{relationhf}
\mathcal{H}=\mathcal{F}\circ\mathcal{A}.
\end{equation}
The next section is devoted to the Kunze-Stein phenomenon on $\Gamma$.  To this end, we will establish the Plancherel formula and the inversion formula of the Helgason-Fourier transform.
Recall that the Plancherel formula and the inversion formula for the spherical Fourier transform were obtained  by different methods in  \cite{fp}, \cite{ks}. More precisely, Iozzi and Picardello \cite{ip1} showed that the Plancherel measure $\mu$ on $\Gamma$ is supported on the set $D\cup E$, where $D$ is the segment $[\frac{\sigma-2Q^{\frac12}}{r(k-1)},\frac{\sigma+2Q^{\frac12}}{r(k-1)}]$ and
$E$ is empty when $k\le r$ and is equal to $\{\frac1{1-k}\}$ if $k>r$. It turns out that the Plancherel measure $\mu$ is the only measure supported on the set $D\cup E$, that verifies the following~:~
$$
f(o)=\int_{D\cup E} \mathcal{H}f(\lambda)\,d\mu(\lambda).
$$
On the other hand, from \cite{fp}, \cite{ks}, we have the following expression of the Plancherel measure
\begin{equation}\label{f(e)}
f(o)= \frac{1}{2\pi}\,\frac{q \ln (q)}{r(k-1)}\,\int_{0}^{\frac{\tau}{2}}\mathcal{H}f(\lambda)\,|c(\lambda)|^{-2}\,d\lambda+\,\!\left[\!\,\frac{(k-r)_{+}}{k}\, \mathcal{H}f(\lambda_0)\right]	
\end{equation}
where the usual notation $(k-r)_+$ stands to 0 if $k\le r$ and to $(k-r)$ if $k>r$, and $\gamma(\lambda_0)=\frac1{1-k}$.  Moreover , $c(\lambda)$ is the meromorphic function given as follows~:~
\begin{equation*}\label{cfunction}
c(\lambda)=\frac{q^{\frac12}}{r(k-1)}\,\frac{q^{\frac12+i\lambda}-(k-1)q^{-\frac12 -i\lambda}+\sigma}{q^{i\lambda}-q^{-i\lambda}}\quad \forall\;\lambda\in \bb C\setminus (\frac{\tau}2)\bb Z.
\end{equation*}
We then can deduce the following~:~
\begin{itemize}
\item
{\bf Plancherel formula}~: For $f\in \mathcal{D}(\mathcal{V})^{\sharp}$, we have
\begin{equation}\label{planch2}
 \left\|f\right\|_{L^2}^{2}=\, \frac{1}{2\pi}\,\frac{q \ln (q)}{r(k-1)}\,\dsps\int_{0}^{\frac{\tau}{2}}|\mathcal{H}f(\lambda)|^2\,|c(\lambda)|^{-2}\,d\lambda\!+\!\,\left[\frac{(k-r)_{+}}{k}\, |\mathcal{H}f(\lambda_0)|^2\right]
\end{equation}
\item
\textbf{Inversion formula}~:\,
For all  $f\in\mathcal{D}(\mathcal{V})^{\sharp}$ and $x\in\mathcal{V}$, we have
\begin{align}\label{inversion}
f(x) &=\, \frac{1}{2\pi}\,\frac{q \ln (q)}{r(k-1)}\,\dsps\int_{0}^{\frac{\tau}{2}}\mathcal{H}f(\lambda)\,\varphi_\lambda(x)\,|c(\lambda)|^{-2}\,d\lambda\\& \!+\!\,\left[\frac{(k-r)_{+}}{k}\, \mathcal{H}f(\lambda_0)\,\varphi_{\lambda_0}(x)\right]\nonumber.
\end{align}
\end{itemize}
We will now generalize in the usual way these results to the nonradial case.  The proof is similar to that of the inversion formula for the Helgason Fourier transform on a symmetric space. More precisely, we first obtain the expression for $f(o)$ using the inversion formula established above for the spherical Fourier transform. Then we derive the expression of $f(x)$ by noting that $f(x)$ is the value at $o$ of a suitable translate of $f$. To this end, we define the spherical means of a function $f:G\to\bb C$ by 
\begin{equation}\label{sphericalmeans}
f^{\sharp}(z)= \frac1{\delta(z)}\sum_{y:\,|y|=|z|}f(y).
\end{equation}
The mean operator $f\to f^{\sharp}$ is the projection over radial functions. It verifies
\begin{equation}\label{expectation}
<f^{\sharp},g^{\sharp}>=\,< f^{\sharp},\,g>=\,< f,g^{\sharp}>,
\end{equation}
$$< f,\, g>=\,\sum_{x\in\mathcal{V}}f(x)\,g(x).$$
Let $f\in\mathcal{D}(\mathcal{V})$. Since $\varphi_\lambda$ is radial, we have
\begin{equation}\label{fphi}
< f^{\sharp},\varphi_\lambda>=\,< f,(\varphi_\lambda)^{\sharp}>=\,< f,\varphi_\lambda>.
\end{equation}
Applying $(\ref{fphi})$ and $(\ref{sphericalmeans})$ to $f^{\sharp}$, we have
\begin{equation}\label{efo}
f(o)=\,f^{\sharp}(o)=\,\frac{1}{2\pi}\,\frac{Q\,\ln(Q)}{r(k-1)}\dsps\int_{0}^{\frac{\tau}{2}}< f,\varphi_\lambda>\,|c(\lambda)|^{-2}\,d\lambda+\,\left[\frac{(k-r)_{+}}{r}\,< f,\varphi_{\lambda_0}>\right].
\end{equation}
Recall the principal series of representations $\pi_{\lambda}$ of G given in \cite{ip1}. For all $\eta\in L^2(\Omega,d\nu)$, 
\begin{equation}\label{representation}
\pi_\lambda(x)\eta(\omega)=\,P(x,\omega)^{\frac12 \!+\!i\lambda}\,\eta(x^{-1}\omega)\hspace{1cm}\forall\;x\in G,\,\forall\;\omega\in\Omega,
\end{equation} 
where $P(x,\omega)$ denotes the Poisson kernel, and $\lambda\in\bb C$. These representations were defined and studied on homogeneous trees in \cite{ftp}, and then on symmetric graphs in \cite{ip1}. Note also that $\pi_{\lambda}$ stands also to the representation on the convolution algebra $\mathcal{D}(\mathcal{V})$as follows,
$$
\pi_{\lambda}(f)=\sum_{x\in\mathcal{V} }f(x)\pi_{\lambda}(x).
$$
Spherical functions are matrix coefficients of $\pi_{\lambda}$. They verify
$$
\varphi_{\lambda}(x)=(\pi_{\lambda}(x)\1,\1)
$$
where $(.|.)$ is the scalar product in $L^2(\Omega,d\nu)$, and $\1$ is the constant function that equals to $1$ on $\Omega$. We now can define the Helgason--Fourier transform in terms of $\pi_{\lambda}$ by
\begin{equation}\label{hftpi}
\hat{f}(\lambda,\omega)=\sum_{x\in\mathcal{V}}f(x)\,P(x,\omega)^{\frac12+i\lambda}=[\pi_{\lambda}(f)\1](\omega).
\end{equation}
\begin{lemma}\label{phietP}
for all $x,y\in G$ we have
$$
\varphi_{\lambda}(x y^{-1})=\varphi_\lambda(y^{-1}x)=\,\int_\Omega P(x,\omega)^{\frac12+i\lambda}\,P(y,\omega)^{\frac12-i\lambda}\,d\nu(\omega).
$$
\end{lemma}
\begin{proof*}
From the cocycle identity $(\ref{cocycle})$, we have, for all $x,y\in G$ and $\omega\in\Omega,$
$$P(y^{-1}x,\omega)=P(x,y\omega)\,P(y^{-1},\omega)\hspace{0.5cm}\text{and}\hspace{0.5cm}1=P(o,\omega)=P(y,y\omega)\,P(y^{-1}\omega),$$
whence
$$
P(y^{-1}x,\omega)=P(x,y\omega)\,P(y,y\omega)^{-1}.
$$
We may then deduce that
\begin{align*}
\varphi_\lambda(y^{-1}x)&=\,\int_\Omega P(x,y\omega)^{\frac12+i\lambda}P(y,y\omega)^{-\frac12-i\lambda}\,d\nu(\omega)\\&=\int_\Omega P(x,y\omega)^{\frac12+i\lambda}\,P(y,y\omega)^{-\frac12-i\lambda}\,\frac{d\nu(\omega)}{d\nu(y\omega)}\,d\nu(y\omega)\\&=\int_\Omega P(x,\omega)^{\frac12+i\lambda}P(y,\omega)^{-\frac12-i\lambda}\,P(y,\omega)\,d\nu(\omega)\\&=\int_\Omega P(x,\omega)^{\frac12+i\lambda}\,P(y,\omega)^{\frac12-i\lambda}\,d\nu(\omega),
\end{align*}
where we have used the fact that $P(y,\omega)=\frac{d\nu(y^{-1}\omega)}{d\nu(\omega)}$ is a Radon-Nikodym derivative.
\end{proof*}
\begin{lemma}\label{lemme}
Let $k\le r$. In this case we have the following~:
\vspace{-12pt}
\begin{itemize}
 \item[(a)]
 \emph{\textbf{Plancherel formula}}~:
For all $f\in\mathcal{D}(\mathcal{V}),$ we have
\begin{equation}\label{planchh}
\dsps\sum_{x\in\mathcal{V}}|f(x)|^2=\frac1{2\pi}\,\frac{q\ln q}{r(k-1)}\int_{0}^{\frac{\tau}2}\int_{\Omega}|\hat{f}(\lambda,\omega)|^2\,|c(\lambda)|^{-2}\,d\nu(\omega)\,d\lambda.
\end{equation}
\item[(b)]
\emph{\textbf{Inversion formula}}~:
For all $f\in\mathcal{D}(\mathcal{V})$ and $x\in\mathcal{V},$ we have
\begin{equation}
 f(x)=\frac1{2\pi}\,\frac{q\ln q}{r(k-1)}\int_{0}^{\frac{\tau}2}\int_{\Omega}\hat{f}(\lambda,\omega)\,P(x,\omega)^{\frac12-i\lambda}\,|c(\lambda)|^{-2}\,d\nu(\omega) d\lambda.
\end{equation}
\end{itemize}
Note that, when $k>r$, there is an additive term that appears on the right hand side of the equality, which corresponds to the parameter $\lambda_0$.
\end{lemma}
\begin{proof*}
The proof is similar to the symmetric space case (see \cite{hel}) and to the homogeneous space case (see \cite{ftp}). More precisely, we apply $(\ref{f(e)})$ on $(f*f^*)^{\sharp}$, where $f^{*}(x)=\overline{f(x^{-1})}$. On one hand we have
$$
(f*f^*)(o)= \dsps\sum_{x\in G}|f(x)|^2=\left\|f\right \|^2_{L ^2(\mathcal{V})}.
$$
On the other hand, from lemma $\ref{phietP}$, we have
\begin{align*}
\mathcal{H}\left[(f*f^*)^{\sharp}\right](\lambda)&=\sum_{x,y\in G}f(x)\overline{f(y)}\varphi_{\lambda}(xy^{-1})\\&=\int_{\Omega}\left[\sum_{x\in G}f(x)\,P(x,\omega)^{\frac12+i\lambda}\right]\times\left[\sum_{y\in G}f(y)\, P(y,\omega)^{\frac12+i\lambda}\right]\,d\nu(\omega)\\&=\int_{\Omega}|\hat{f}(\lambda,\omega)|^2\,d\nu(\omega)
\end{align*}
for all $\lambda\in\bb R$. We conclude by using $(\ref{phietP})$, proven (a). In order to prove (b), we apply $(\ref{f(e)})$ to $f_x^{\sharp}$. On one hand, $(f_x)^{\sharp}(o)=f(x)$. On the other hand, using $(\ref{sphericalmeans})$ and lemma $\ref{phietP}$, one can show that
\begin{align*}
\mathcal{H}\left[(f_x)^{\sharp}\right](\lambda)&=\sum_{\dsps\mathop{x,y\in G}_{|y|=|x|}}f(xy)\frac{\varphi_{\lambda}(x,y)}{\delta(|z|)}\\&=\sum_{y\in G}f(y)\varphi_{\lambda}(x^{-1}y)=\int_{\Omega}\left[\sum_{y\in G}f(y)P(y,\omega)^{\frac12+i\lambda}\right]P(x,\omega)^{\frac12-i\lambda}\\&=\int_{\Omega}\hat{f}(\lambda,\omega)\,P(x,\omega)^{\frac12-i\lambda}d\nu(\omega)
\end{align*}
for all $\lambda\in\bb R$. We conclude using $(\ref{f(e)})$.
\end{proof*}
Now let $f\in\mathcal{D}(\mathcal{V})$ and $\chi\in\mathcal{D}(\mathcal{V})^{\sharp}$. One can easily show that
\begin{equation}\label{convolution1}
\widehat {(f*\chi)}(\lambda,\omega)=\,\hat{f}(\lambda,\omega)\,\mathcal{H}\chi(\lambda).
\end{equation}
Using the Plancherel formula $(\ref{planchh})$, and the fact that $|\varphi_{\lambda}(x)|\le |\varphi_{0}(x)|\;\forall\,\lambda\in\bb R,\,\forall\,x\in\mathcal{V},$ one can show that
$$
\left\|f*\chi\right\|_{L^2}\le \left\|f\right\|_{L^2}\,\sum_{x\in\mathcal{V}}\chi(x)\varphi_0(x).
$$
It is easily shown, using $(\ref{phi})$, that for all $x\in\Gamma$ we have $\varphi_0(x)\le C\,(1+|x|)\,Q^{-\frac{|x|}2}\,\forall\,x\in\mathcal{V}$. Then, the following version of the Kunze--Stein phenomenon holds~:
\begin{equation}\label{ks}
\left\|f*\chi\right\|_{L^2(\mathcal{V})}\le C\,\left\|f\right\|_{L^2(\mathcal{V})}\,\dsps\sum_{n\ge 0}\chi(n)\,(1\!+\!n)\,Q^{\frac{n}{2}}.
\end{equation}
In the next proposition we  generalize this result to general $L^p$ spaces.
\begin{proposition}
For all $2\le p,\tilde{p}<\infty,$ there exist $C>0$ such that
$$
\left\|f*\chi\right\|_{L^p}\le C\,\left\|f\right\|_{L^{{\tilde{p}}^{\prime}}}\,\Big\{\dsps\sum_{n\ge 0}|\chi(n)|^{r}	\,(1+n)^{2s}\,q^{(1+s)n})\Big\}^{\frac{1}{r}},
$$
where $r=\frac{p\,\tilde{p}}{p+\,\tilde{p}}$ and $s=\,\frac{\operatorname{min}\,\left\{p,\tilde{p}\right\}}{p+\,\tilde{p}}.$ 
\end{proposition}
The same result has been obtained on symmetric spaces and more generally on Damek-Ricci spaces in \cite{ady}. The proof is based on an interpolation argument between the following inequalities,
$$
\left\|f*\chi\right\|_{L^p}\le\,\left\|f\right\|_{L^1}\left\|\chi\right\|_{L^p},\;\;\;\left\|f*\chi\right\|_{L^{\infty}}\le\,\left\|f\right\|_{L^{{\tilde{p}}^{\prime}}}\left\|\chi\right\|_{L^p}.
$$
\noindent
The interpolation argument is well detailed in (\cite{ady}, lemma 5.1).

\section{The shifted wave equation}
 This section is devoted to the shifted wave equation on a symmetric graph $\Gamma$. We will first introduce some notations. 
 Let $\alpha, \,\beta$ be two reals define as follows~:~
 $$
 \alpha=\frac{Q+1}{r(k-1)}\;\text{and}\;\beta=\frac{2Q^{\frac12}}{r(k-1)}.
 $$
 Remark that $\alpha-\beta=1-\gamma(0)$ is the $L^2$ spectral gap of the Laplacian. 
 On radial functions, the Laplacian is given by 
 $$\mathcal{L}^{\Gamma}f(0)=\,f(0)-f(1)$$
and
\begin{equation}\label{laplacien radial}
\mathcal{L}^{\Gamma}f(n)=\frac{1}{r(k-1)}\{(Q+1)f(n)-f(n-1)-Q\,f(n+1)\}
\end{equation}
for $n\in \N^*$. On horocyclical functions, the Laplacian is given by 
\begin{align}\label{laplacien horocyclique}
\mathcal{L}^{\Gamma}f(h)&=\,\frac{1}{r(k-1)}\{(Q+1)f(h)-Q\,f(h-1)-f(h+1)\}\\&=\,\beta\,Q^{\frac{h}{2}}\mathcal{L}^{\bb Z}_h\{Q^{-\frac{h}{2}}\,f(h)\}+(1-\gamma(0))f(h)\nonumber\\&=\,\beta\,Q^{\frac{h}{2}}\mathcal{L}^{\bb Z}_h\{
Q^{-\frac{h}{2}}\,f(h)\}+(\alpha-\beta)f(h).\nonumber
\end{align}
 We are interested in the following wave equation on $\Gamma$~:
\begin{equation}\label{wave}
\begin{cases}
\beta\,\mathcal{L}^{\bb Z}_n\,u(x,n)=\,(\mathcal{L}_x^{\Gamma}-(\alpha-\beta))\,u(x,n),\\
u(x,0)=\,f(x),\;\{u(x,1)-u(x,-1)\}/2=\,g(x).
\end{cases}
\end{equation}
We will solve this equation by using the inverse dual Abel transform and the following discrete version of the \`Asgeirsson theorem.
\begin{theorem}\label{asgeirsson}
Let $U$ be a function on $\Gamma\times\Gamma$ such that 
\begin{equation}\label{U}
\mathcal{L}^{\Gamma}_x U(x,y)=\mathcal{L}^{\Gamma}_y U(x,y)\;\;\forall\,x,y\in\Gamma.
\end{equation}
Then
\begin{equation}\label{moyenneA}
\sum_{x^{'}\in S(x,m)}\dsps\sum_{y^{'}\in S(y,n)}U(x^{'},y^{'})=\sum_{x^{'}\in S(x,n)}\sum_{y^{'}\in S(y,m)}U(x^{'},y^{'})
\end{equation}
for all $x,y\in\Gamma$ and $m,n\in \bb N$. In particular
\begin{equation}\label{moyenneA1}
\sum_{x^{'}\in S(x,n)}U(x^{'},y)=\sum_{y^{'}\in S(y,n)}U(x,y^{'}).
\end{equation}
\end{theorem}
In order to prove Theorem $\ref{asgeirsson}$ we will need the following lemma.
\begin{lemma}\label{lemmee}
Consider the spherical means 
$$
f^{\sharp}_x(n)=\frac{1}{\delta(n)}\sum_{y\in S(x,n)}f(y)\;\;\;\forall x\in\Gamma,\,\forall n\in\bb N.
$$
Then 
$$
(\mathcal{L}^{\Gamma}f)^{\sharp}_x(n)=\,(\operatorname{rad}\,\mathcal{L}^{\Gamma})_n\,f^{\sharp}_x(n),
$$
where $\operatorname{rad}\,\mathcal{L}$ denotes the radial part $(\ref{laplacien radial} )$ of $\L^{\Gamma}$.
\end{lemma}
\begin{proof*}[Proof of Lemma \emph{\ref{lemmee}}]
We have 
\begin{equation*}
(\mathcal{L}^{\Gamma}f)^{\sharp}_x(n)=\begin{cases}
f(x)-f^{\sharp}_x(1)\hspace{5.6cm}\qquad\textrm{if}\,n=0,\\
\frac{1}{r(k-1)}\{(Q+1)f^{\sharp}_x(n)-f^{\sharp}_x(n-1)-Q\,f^{\sharp}_x(n+1)\}\quad\textrm{if}\,n\in\bb{N}^*.
\end{cases}
\end{equation*}
\end{proof*}

\begin{proof*}[Proof of Theorem \emph{\ref{asgeirsson}}]
Let $x,y\in\Gamma$. We define the double spherical means of $U$
\begin{equation*}\label{ds}
U^{\sharp\sharp}_{x,y}(m,n)=\frac{1}{\delta(m)}\sum_{x^{'}\in S(x,m)}\frac{1}{\delta(n)}\sum_{y^{'}\in S(y,n)}U(x,y),
\end{equation*}
and we denote it by $V(m,n)$. We have
\begin{equation}\label{V}
(\operatorname{rad} \mathcal{L}^{\Gamma})_m V(m,n)=\,(\operatorname{rad} \mathcal{L}^{\Gamma})_n V(m,n).
\end{equation}
Let 
$$M^{(n)} f(x)=\,f^{\sharp}_x(n)$$
denotes the mean operator, then 
$$U^{\sharp\sharp}_{x,y}(m,n)=\,M^{(m)}_{x}\,M^{(n)}_{y} U(x,y).$$
According to  $(\ref{U})$ and lemma $\ref{lemmee}$, we have
\begin{align*}
(\operatorname{rad}\mathcal{L}^{\Gamma})_{m} \,M^{(m)}_{x}\,M^{(n)}_{y}\,U(x,y)&=\,M^{(m)}_x\,\mathcal{L}^{\Gamma}_x\,M^{(n)}_y U(x,y)=\,M^{(m)}_x\,M^{(n)}_y\,\mathcal{L}^{\Gamma}_x U(x,y)\\&=M^{(m)}_x\,M^{(n)}_y\,\mathcal{L}^{\Gamma}_y\,U(x,y)=(\operatorname{rad}\mathcal{L}^{\Gamma})_n\,M^{(m)}_x\,M^{(n)}_y U(x,y),
\end{align*}
thus $(\ref{V})$ holds.
Let us prove the symmetry
\begin{equation}\label{symetry}
V(m,n)=\,V(n,m)\;\;\;\forall\,m,n\in\bb N
\end{equation}
by reccurence on  $\ell=m+n.$ First of all, $(\ref{symetry})$ is trivial if $\ell=0$, and for $\ell=1$, $(\ref{symetry})$ results from the definition of $\operatorname{rad} \,\mathcal{L}^{\Gamma}$ and from $(\ref{U})$. Assume next $\ell\ge 1$ and that $(\ref{symetry})$ holds for $m+n\le \ell$.  Let $m>n>0$, with $m+n=\ell+1$ and let $1\le \ell^{\prime}\le m-n$. We deduce from $(\ref{V})$ at the point  $(m-\ell^{\prime},n+\ell^{\prime}-1)$ that 
\begin{align}
&V(m-\ell^{\prime}+1,n+\ell^{\prime}-1)-V(m-\ell^{\prime},n+\ell^{\prime})\label{mkn}\\&=\frac{1}{Q}\,\{V(m-\ell^{\prime},n+\ell^{\prime}-2)-V(m-\ell^{\prime}-1,n+k-1)\},\nonumber
\end{align}
Adding up $(\ref{mkn})$ over $\ell^{\prime}$, we obtain
\begin{equation}\label{mkn1}
V(m,n)-V(n,m)=\,\frac{1}{Q}\,\{V(m-1,n-1)-V(n-1,m-1)\}.
\end{equation}
which vanishes by induction. Now, using $(\ref{V})$ at the points $(0,\ell)$ and $(\ell, 0)$, we deduce that
\begin{equation*}
\begin{cases}
V(0,\ell+1)=\,\tfrac1Q\,\{-(k-2)V(0,\ell)-V(0,\ell-1)+r\,(k-1)V(1,\ell)\},\\
V(\ell+1,0)=\,\tfrac1Q\,\{-(k-2)V(\ell,0)-V(\ell-1,0)+r\,(k-1)V(\ell,1)\}.
\end{cases}
\end{equation*}
Hence $V(\ell+1,0)=V(0,\ell+1)$ by using $(\ref{mkn1})$ and by induction. This concludes the proof of Theorem $\ref{asgeirsson}$.
\end{proof*} 
We are now ready to solve explicitely the shifted wave equation $(\ref{wave})$ on $\Gamma$. Let us first consider a solution $u$ to $(\ref{wave})$ with initial data 
$$
u(x,0)=f(x),\;\{u(x,1)-u(x,-1)\}/2=0.$$
First of all, we have that the function $(x,n)\mapsto u(x,-n)$ satisfies the same Cauchy problem, then $u(x,n)=u(x,-n)$  by uniqueness. Now, according to $(\ref{laplacien horocyclique})$, the function $$U(x,y)=\,Q^{\frac{h(y)}{2}}\,u(x,h(y))\;\;\forall\,x,y\in\Gamma$$
satisfies $(\ref{U})$. By applying $(\ref{moyenneA})$ to $U(x,y)$ at $y=0$, we deduce that the dual Abel transform of $n\mapsto u(x,n) $ is equal to the spherical means $f^{\sharp}_{x}(n)$ of the initial data $f$. Thus 
 \begin{equation*}\label{solution1}
 u(x,n)=\,(\mathcal{A}^*)^{-1}(f^{\sharp}_x)(n)\;\;\forall\,x\in\Gamma,\,\forall\,n\in\bb N.
 \end{equation*}
Consider next a solution to the Cauchy problem $(\ref{wave})$ with initial data
\begin{equation*}
u(x,0)=0 \;\;\textrm{and}\;\; \!\frac{u(x,1)-u(x,-1)}{2}=g(x).
\end{equation*} 
   Then $n\mapsto u(x,n)$ is an odd function, and 
   $$
   v(x,n)=\frac{u(x,n+1)-u(x,n-1)}{2}
   $$ 
   is a solution to $(\ref{wave})$ with initial data
  \begin{equation*}
v(x,0)=g(x)\;\;\textrm{and}\;\; \frac{v(x,1)-v(x,-1)}{2}=0.
\end{equation*} 
Hence 
\begin{equation*}
u(x,n)=\begin{cases}
2\,\operatorname{sign}(n)\,\sum_{0 < \ell\,\textrm{odd} < |n|}v(x,\ell)\hspace{0.6cm}\qquad &\textrm{if}\,n\,\textrm{is even},\\
g(x)+\,2\,\operatorname{sign}(n)\,\sum_{0< \ell\,\textrm{even}<|n|}v(x,\ell)&\textrm{if}\,n\,\textrm{is odd},
\end{cases}
\end{equation*}
with $v(x,\ell)=\,(\mathcal{A}^*)^{-1}(g^{\sharp}_x)(\ell).$ Thus the solution to the Cauchy problem $(\ref{wave})$ is given by 
\begin{equation*}
u(x,n)=
(\mathcal{A}^{*})^{-1}(f^{\sharp}_x)(n)+\,2\,\operatorname{sign}(n)\sum_{0< \ell\,\textrm{odd}< |n|}(\mathcal{A}^{*})^{-1}(g^{\sharp}_x)(\ell)
\end{equation*}
if\,$n$\,is even, and
\begin{equation*}
u(x,n)=(\mathcal{A}^{*})^{-1}(f^{\sharp}_x)(n)+\,g(x)+\,2\,\operatorname{sign}(n)\,\sum_{0< \ell\,\textrm{even}< n}(\mathcal{A}^{*})^{-1}(g^{\sharp}_x)
\end{equation*}
if $n$ is odd. Using the inverse dual Abel transform, we are now able to give an explicite expression of the solution to $(\ref{wave})$.
\vspace{-0.4cm}
\vskip12pt
\begin{theorem}
\underline{When $k<r$},  the solution to the Cauchy problem $(\ref{wave})$ is given by\\
$u(x,0)=f(x)$, and 
\begin{align}\label{solutionwave}
&u(x,n)=\\&\,\frac12\,Q^{-\frac{|n|}2}\,\sum_{d(x,y)=|n|}f(y)\nonumber\\&-\,\frac1{2k}\,Q^{-\frac{|n|}2}\sum_{0\le\ell<|n|}\{Q-1+(r-k)(1-k)^{|n|-\ell}\}\sum_{d(x,y)=\ell}f(y)\nonumber\\&+\,\operatorname{signe}(n)\,Q^{-\frac{|n|-1}{2}}\sum_{d(x,y)=|n|-1}g(y)\nonumber\\&+\,\operatorname{signe}(n)\,\frac{1}{k}\,Q^{-\frac{|n|-1}{2}}\,\left\{\sum_{0\le d(x,y)<|n|-1}g(y)-\,\sum_{0\le\ell< |n|-1}(1-k)^{|n|-\ell}\sum_{d(x,y)=\ell}g(y)\right\}\nonumber
\end{align}
if $n\in\bb Z^*.$\\
{\underline{When $k=r$}},  the solution to $(\ref{wave})$ is given by\\
 $u(x,0)=f(x)$ and
\begin{align}\label{solutionk=r}
 &u(x,n)=\frac12\,(k-1)^{-|n|}\,\sum_{d(x,y)=|n|}f(y)\\&-\frac{k-2}2\,(k-1)^{-|n|}\,\sum_{d(x,y)<|n|}f(y)+\,\operatorname{signe}(n)\,(k-1)^{-(|n|-1)}\,\sum_{d(x,y)=|n|-1}g(y)\nonumber\\&+\operatorname{signe}(n)\,\frac1{k}\,(k-1)^{-(|n|-1)}\times\nonumber\\&\times\left\{\sum_{d(x,y)<|n|-1}g(y)+\sum_{0\le\ell<|n|-1}(k-1)^{|n|-\ell}\sum_{d(x,y)=\ell}g(y)\right\}\nonumber
\end{align}
if $n\in\bb Z^*.$
\begin{proof*}
The case $n=0$ is trivial. Let $n\in\N^*$ and assume that $k<r$. For $f~:~\Gamma\to\C$ and $j\in\N^*$, we introduce the functions 
$$
f_j(x)=\sum_{d(x,y)=j}f(y).
$$
Using the expression of the inverse dual Abel transform, we have 
\begin{align*}
&(\mathcal{A}^*)^{-1}(f^{\sharp}_x)(n)=\frac12\,Q^{-\frac{n}2}\,f_n(x)\\&-\frac1{2k}\,Q^{-\frac{n}2}\,\sum_{0\le j< n}\{Q-1+(r-k)(1-k)^{n-j}\}f_j(x).\nonumber
\end{align*}
When $n$ is even, we have
\begin{align*}
2\sum_{0<\ell\;\rm{odd}\,< n}(\mathcal{A}^*)^{-1}(g_x^{\sharp})(\ell)&=\sum_{0<\ell\;\rm{odd}\,< n}Q^{-\frac{\ell}2}\,g_{\ell}(x)\\&-\frac1{k}\,\sum_{0<\ell\;\rm{odd}\,< n}Q^{-\frac{\ell}2}\sum_{0\le j<\ell}\{Q-1+(r-k)(1-k)^{\ell-j}\}\,g_{\ell}(x).
\end{align*}
\vspace{-0.3cm}
We decompose this expression, by separating it with respect to the parity of $j$ 
as follows~:
\begin{align}
&2\sum_{0<\ell\;\rm{odd}\,< n}(\mathcal{A}^*)^{-1}g_x^{\sharp}(\ell)\nonumber =Q^{-\frac{n-1}2}\sum_{d(x,y)=n-1}g(y)+\sum_{0<\ell\;\rm{odd}\,< n-1}Q^{-\frac{\ell}2}\,g_{\ell}(x)\\&-\frac{Q-1}{k}\sum_{0<\ell\;\rm{odd}\,< n}Q^{-\frac{\ell}2}\sum_{0< j\;\rm{odd}\,< \ell}\,g_j(x)\label{20}\\&-\frac{Q-1}{k}\sum_{0<\ell\;\rm{odd}\,< n}Q^{-\frac{\ell}2}\sum_{0\le j\;\rm{even}<\ell}\,g_j(x)\label{21}\\&-\frac{r-k}{k}\sum_{0<\ell\;\rm{odd}\,< n}(1-k)^{\ell}\,Q^{-\frac{\ell}2}\sum_{0< j\;\rm{odd}\,<\ell}(1-k)^{-j}\,g_j(x)\label{22}\\&-\frac{r-k}{k}\sum_{0<\ell\;\rm{odd}\,< n}(1-k)^{\ell}\,Q^{-\frac{\ell}2}\sum_{0\le j\;\rm{even}\,<\ell}(1-k)^{-j}\,g_j(x)\label{23}
\end{align}
We may simplify this expression. On one hand
\begin{align}
(\ref{20})&= -\frac{Q-1}{k}\sum_{0< j\;\textrm{ even}< n-1}\hspace{0.3cm}\sum_{j<\ell\;\textrm{ even}\,< n}Q^{-\frac{\ell}2}\,g_j(x)\nonumber\\&=-\frac1{k}\sum_{0< j\;\rm{ even}< n-1}Q^{-\frac{j}2}\,g_j(x)+\frac{Q^{-\frac{n-1}2}}{k}\sum_{0<j\;\rm{ even}<n-1}g_j(x)\nonumber
\end{align}
by using the geometric sum
$$\sum_{j<\ell\;\rm{even}<n}Q^{-\frac{\ell}2}=\,\frac{Q^{-\frac{j}2}-Q^{-\frac{n-1}2}}{Q-1}.$$
Similarly
\begin{equation*}
(\ref{21})=-\frac1{k}\sum_{0\le j\;\rm{even}\,< n}Q^{-\frac{j-1}2}\,g_j(x)+\frac{Q^{-\frac{n-1}{2}}}{k}\sum_{0\le j\;\rm{even}<n}g_j(x).
\end{equation*}
\noindent\ssb
On the other hand 
\begin{align*}
(\ref{22})&=-\frac{r-k}{k}\sum_{0< j\;\rm{odd}\,<n-1}(1-k)^{-j}\sum_{j< \ell\;\rm{odd}< n}(1-k)^{\ell}\,Q^{-\frac{\ell}2}\,g_j(x)\\&=\frac{1}{k}\sum_{0< j\;\rm{odd}\,<n-1}\left\{(1-k) Q^{-\frac{j}2}-(1-k)^{n}Q^{-\frac{n-1}2}\,g_j(x)\right\}\\&=-(1-\frac1{k})\sum_{0< j\;\rm{odd}\,< n-1}Q^{-\frac{j}{2}}\,g_j(x)\\&-\frac{Q^{-\frac{n-1}2}}{k}\sum_{0< j\;\rm{odd}\,< n-1}(1-k)^{n-j}\,g_j(x)
\end{align*}
by using the geometric sum
$$\sum_{j<\ell\;\rm{odd}\,<n}\left[\frac{(1-k)^2}{Q}\right]^{\frac{\ell}2}=\,-\frac1{r-k}\left\{(1-k)^{j+1}\,Q^{-\frac{j}2}-(1-k)^n\,Q^{-\frac{n-1}2}\right\}.$$
\noindent
Similarly, we have 
\begin{align*}
(\ref{23})&=-\frac{r-k}{k}\sum_{0\le j\;\rm{even}\,< n} (1-k)^{-j}\sum_{j+1\le\ell\le n-1}(1-k)^{\ell}Q^{-\frac{\ell}2}\,g_j(x)\\&=\frac{1}{k}\sum_{0\le j\;\rm{even}\,< n}\left\{Q^{-\frac{j-1}2}-(1-k)^{n-j}Q^{-\frac{n-1}2}\right\}\,g_j(x)\\&=\frac1{k}\sum_{0\le j\;\rm{even}\,< n}Q^{-\frac{j-1}2}\sum_{d(x,y)=j}g(y)-\frac{Q^{-\frac{n-1}2}}{k}\sum_{0\le j\;\rm{even}\,< n}(1-k)^{n-j}\,g_j(x).
\end{align*}
When $n$ is odd, a similar computation gives the same expression. Hence $(\ref{solutionwave})$ holds for $k<r$. \\
When $k=r$,  $(\ref{22})$ and $(\ref{23})$ vanishes, and we get the expression $(\ref{solutionk=r})$.

\end{proof*}

\end{theorem}


%


\begin{thebibliography}{99}

\bibitem{ady}
J.-Ph. Anker, E. Damek, C. Yacoub~:
\textit{Spherical analysis on harmonic $AN$ groups\/},
Ann. Scuola Norm. Sup. Pisa 23 (1996), 643--679.
\bibitem{amps}
J.-Ph. Anker, P. Martinot, E. Pedon, A.G. Setti~:
\textit{The shifted wave equation on Damek-Ricci spaces and homogeneous trees\/},
dans \textit{Trends in harmonic analysis}, Springer INdAM Ser. 3 (2013), 1--25.
\bibitem{bfp}
W. Betori, J. Faraut, M. Pagliacci~:
\textit{An inversion formula for the Radon transform on trees\/},
Math. Z. 201 (1989), 327--337.
\bibitem{bp}
W. Betori, M. Pagliacci~:
\textit{Harmonic analysis for groups acting on trees\/},
Boll. Unione Mat. Ital. 6 (1984), 333--345 
\bibitem{bp1}
W. Betori, M. Pagliacci~:
\textit{The Radon transform on trees\/},
Boll. Un. Mat. Ital. B (6), vol. 5, (1986), 267--277.

\bibitem{ca}
P. Cartier~: \textit{G\'eom\'etrie et analyse sur les arbres\/},
 S\'eminaire Bourbaki, 1971--1972, expos\'e 407, Lect. Notes Math. 317, Springer (1973), 123--140.
\bibitem{ca1}
P. Cartier~: 
\textit{Harmonic analysis on trees\/},
Proc. Symp. Pure Math. A.M.S. 26 (1972), 419--424 
\bibitem{cas}
D.I. Cartwright, P.M. Soardi~:
\textit{Harmonic analysis on the free product of two cyclic groups\/},
 J. Funct. Anal. 65 (1986), no. 2, 147--171.
 \bibitem{ccc}
 E. Casadio Tarabusi, J.M. Cohen, et F. Colonna~:
\textit{The horocyclic Radon transform on nonhomogeneous trees\/},
 Israel J. Math. 78 (1992), no. 2-3, 363--380.
 \bibitem{cms}
 M.G. Cowling, S. Meda, A.G. Setti~: \textit{An overview of harmonic analysis on the group of isometries of a homogeneous tree\/}, Expo. Math. 16 (1998), 385--423.
\bibitem{fp}
J. Faraut, M. Picardello~:
\textit{The Plancherel measure for symmetric graphs\/},
 Ann. Mat. Pura Appl. (4) 138 (1984), 151--155.
\bibitem{ftn}
A. Fig\`a-Talamanca, C. Nebbia~: \textit{Harmonic analysis and representation theory for groups acting on homogeneous trees\/}, London Math. Soc. Lect. Notes Ser. 162, Cambridge Univ. Press, 1991.
\bibitem{ftp}
A. Fig\`a-Talamanca, M. Picardello~:
\textit{Spherical functions and harmonic analysis on free groups\/},
 J. Funct. Anal. 47 (1982), no. 3, 281--304.
 \bibitem{hel1}
 S. Helgason~:
 \textit{The {R}adon transform\/},
 Progress in Mat. (5), Birkh\"auser Boston,  1980.
 \bibitem{hel}
S. Helgason~:
\textit{Groups and Geometic Analysis\/},
Academic Press, 1984.
 \bibitem{ip}
A. Iozzi, M. Picardello~:
\textit{Graphs and convolution operators\/},
 Topics in modern harmonic analysis, Vol. I, II (Turin/Milan, 1982), 187--208, Ist. Naz. Alta Mat. Francesco Severi, Rome, 1983.
\bibitem{ip1}
A. Iozzi, M. Picardello~:
\textit{Spherical functions on symmetric graphs\/},
Harmonic analysis (Cortona, 1982), 344--386, Lecture Notes in Math., 992, Springer, Berlin, 1983.
\bibitem{alje}
A. Jamal Eddine~:
\textit{Equations d'\'evolution sur certains groupes hyperboliques\/},
Th\`ese de doctorat, 2013.
\bibitem{ks}
G; Kuhn, P.M. Soardi~:
\textit{The Plancherel measure for polygonal graphs\/},
Ann. Mat. Pura Appl. 4 (1983), 393--401. 
\bibitem{wo1}
W. Woess~:
\textit{Random walks on infinite graphs and groups\/},
Cambridge Tracts in Mathematics 138, Cambridge University Press.







\end{thebibliography}
\end{document}